\documentclass[a4paper]{amsart}
\usepackage[utf8]{inputenc}
\usepackage{amssymb}
\usepackage{amsmath}
\usepackage{euscript}
\usepackage{amscd}
\usepackage{amsthm}
\usepackage[matrix,arrow,curve]{xy}
\usepackage{graphicx}
\DeclareGraphicsExtensions{.jpg}

\newcommand*{\hm}[1]{#1\nobreak\discretionary{} {\hbox{$\mathsurround=0pt #1$}}{}}

\newenvironment{kkeywords}
{
\list{}{\itemindent\listparindent\leftmargin 2\parindent
\rightmargin\leftmargin}
\leftmargini 2em
\small\item[]\relax\parindent12pt\leftmarginii3em\hskip\parindent{Keywords:}}
{\endlist\vskip12pt plus 3pt minus 2pt}

\theoremstyle{problems}

\newtheorem{theorem}{Theorem}[section]
\newtheorem{proposition}[theorem]{Proposition}
\newtheorem{lemma}[theorem]{Lemma}
\newtheorem{corollary}[theorem]{Corollary}

\theoremstyle{definition}
\newtheorem{example}[theorem]{Example}

\begin{document}

\title[Massey products for Pogorelov polytopes]{MASSEY PRODUCTS IN COHOMOLOGY OF MOMENT-ANGLE MANIFOLDS CORRESPONDING TO POGORELOV POLYTOPES}

\author{Elizaveta Zhuravleva}

\address{Department of Mathematics and Mechanics, Moscow
State University, Leninskie Gory, 119991 Moscow, Russia} 
\email{ahertip@gmail.com}

\thanks{The work was supported by the Russian Foundation for Basic Research (grant
no. 18-51-50005).}

\begin{abstract}
In this work we construct nontrivial Massey products in the cohomology of moment-angle manifolds corresponding to polytopes from the Pogorelov class.
This class includes the dodecahedron and all fullerenes, i. e. simple 3-polytopes with only 5-gonal and 6-gonal facets.
The existence of a nontrivial Massey product implies the nonformality of the space in the sense of rational homotopy theory.
\end{abstract}

\maketitle

\vspace{-3ex}

\begin{kkeywords}
\textit{Massey products, moment-angle complexes, Pogorelov polytopes, fullerenes.}
\end{kkeywords}

\bigskip

\section{Introduction}

\medskip

In this paper we consider the problem of existence of nontrivial Massey products in the cohomology of moment-angle manifolds corresponding to 3-dimensional simple polytopes $P$.
As shown in [1], these manifolds $\mathcal{Z}_P$ are 2-connected smooth manifolds of dimension $m+3$, where $m$ is the number of facets in $P$.
First examples of moment-angle manifolds having nontrivial Massey products were found by Baskakov (see [2]).
Limonchenko constructed in [3] a family of moment-angle manifolds having nontrivial Massey $n$-products for any $n$.

In dimension $3$ the Pogorelov class of simple polytopes is particularly interesting.
This class consists of combinatorial 3-dimensional simple polytopes  which do not have 3-belts and 4-belts of facets.
It is known that the Pogorelov class consists precisely of combinatorial 3-polytopes which admit a right-angled realization in Lobachevsky space $\mathbb{L}^3$, and such a realization is unique up to isometry (see [4], [5], [6]).
There is a family of hyperbolic $3$-manifolds associated with Pogorelov polytopes, known as hyperbolic manifolds of L\" obell type (see [7]).
Moment-angle manifolds corresponding to Pogorelov polytopes are important for the topological study of hyperbolic manifolds of L\"obell type, and also for cohomological  rigidity of $6$-dimensional  (quasi)toric  manifolds.

It is known that there are no triple Massey products of 3-dimensional classes in cohomology of moment-angle manifolds corresponding to Pogorelov polytopes (see [6]). For cohomology classes of higher dimension, the existence of Massey products was open.
We prove that for any Pogorelov polytope $P$ the corresponding moment-angle manifold $\mathcal{Z}_P$ has a nontrivial triple Massey product in the cohomology.
This implies that all such manifolds $\mathcal{Z}_P$ are non-formal.

Our construction of nontrivial Massey products is based on the combinatorial description of the cohomology of moment-angle complexes and certain combinatorial properties of Pogorelov polytopes.
The Pogorelov class contains all fullerenes (simple polytopes with only 5-gonal and 6-gonal facets), in particular, the dodecahedron. In this paper we also consider a particular case of dodecahedron.

The author is grateful to the advisor Taras Panov for the formulation of the problem and for his permanent attention to this work.

\bigskip

\section{Preliminaries}

\medskip

Let $A=\bigoplus \limits _{i \geqslant 0} A^i$ be a commutative differential graded algebra over $\mathbb{Z}$.
Let \linebreak $\alpha_i \hm \in H^{k_i}(A), i=1, 2, 3$, be three cohomology classes such that $\alpha_1 \alpha_2 = 0,$ $\alpha_2 \alpha_3 \hm = 0$ in $H(A)$.
Choose their representing cocycles $a_i \in A^{k_i}, \, i=1, 2, 3$.
Since the pairwise cohomology products vanish, there are elements $a_{12} \in A^{k_1+k_2-1}$ and $a_{23} \in A^{k_2+k_3-1}$ such that
$d a_{12} = a_1 a_2 \: \mbox{ and } d a_{23} = a_2 a_3.$
Then one easily checks that
$$b = (-1)^{k_1+1} a_1 a_{23} + a_{12} a_3$$
is a cocycle in $A^{k_1+k_2+k_3-1}$.

A \textit{triple Massey product} $\langle \alpha_1, \, \alpha_2, \, \alpha_3 \rangle$ is the set in   $H^{k_1+k_2+k_3-1}(A)$ consisting of all elements obtained by this procedure.
Since elements $a_{12}$ and $a_{23}$ are defined up to addition of cocycles in $A^{k_1+k_2-1}$ and $A^{k_2+k_3-1}$ respectively, then, more precisely, we have
$$\langle\alpha_1, \alpha_2, \alpha_3\rangle = [b] +\alpha_1 H^{k_2+k_3-1}+\alpha_3 H^{k_1+k_2-1}.$$

The subset $\alpha_1 H^{k_2+k_3-1}+\alpha_3 H^{k_1+k_2-1}$ is called the \textit{indeterminacy} of a Massey product $\langle\alpha_1, \alpha_2, \alpha_3\rangle$.

A Massey product $\langle \alpha_1, \, \alpha_2, \, \alpha_3 \rangle$ is called \textit{trivial} if $0 \in \langle \alpha_1, \, \alpha_2, \, \alpha_3 \rangle$ and \textit{nontrivial} otherwise.

Let $\mathcal{K}$ be a simplicial complex on the set $[m]= \{1, \ldots, m\}$.
The \textit{moment-angle complex} (see [1]) corresponding to a simplicial complex  $\mathcal{K}$ is a topological space defined by
$$\mathcal{Z}_{\mathcal{K}} = \bigcup_{I \in \mathcal{K}} (\prod_{i \in I} D^2) \times ( \prod_{i \notin I} S^1)  \subseteq (D^2)^m.$$

An important class of simplicial complexes  $\mathcal{K}$ comes from simple polytopes.
Recall that an $n$-dimensional polytope $P$ is called \textit{simple} if exactly $n$ facets meet at each vertex.
Denote by $\mathcal{K}_P$ the simplicial complex dual to the boundary of a simple polytope $P$.
In more detail, if $\{F_1, \ldots, F_m\}$ is the set of faces of codimension 1 in $P$, then
$$\mathcal{K}_P=\{\{i_1, \ldots, i_k\} \subseteq [m] : F_{i_1} \cap \cdots \cap F_{i_k} \neq \varnothing \}.$$
Note that $\mathcal{K}_P$ is a triangulation of the $(n-1)$-dimensional sphere.

\medskip

\begin{theorem}[{\cite[Theorem 6.2.4, Corollary 6.2.5]{1}}]
$\mathcal{Z}_{\mathcal{K}}$ is a CW-complex. If  $\mathcal{K}=\mathcal{K}_P$ for a simple $n$-polytope $P$, then
$\mathcal{Z}_{\mathcal{K_P}}$
 is a smooth manifold.
\end{theorem}
\medskip

Let $I=\{i_1, \ldots , i_k\}$ be a simplex, $i_1 < i_2 < \cdots < i_k$.
Denote by $v_I$ the monomials $v_{i_1} \cdots v_{i_k}$ in the polynomial algebra $\mathbb{Z}[v_1, \ldots, v_m]$, and denote by $u_I$ the exterior monomial $u_{i_1} \cdots u_{i_k}$ in the exterior algebra $\Lambda[u_1, \ldots, u_m]$.
The \textit {face ring} of the simplicial complex $\mathcal{K}$ on the set $[m]$ is defined as the
quotient of the polynomial algebra by the monomial ideal corresponding to non-simplices of $\mathcal{K}$:
$$\mathbb{Z}[\mathcal{K}] \hm =\mathbb{Z}[v_1, \ldots, v_m]/ \mathcal{I_K},$$
where $\mathcal{I_K}=(v_I: I\notin \mathcal{K})$ is the Stanley-Reisner ideal.

We define the quotient algebra
$$R^{*}(\mathcal{K})=\Lambda [u_1, \ldots, u_m]\otimes \mathbb{Z}[\mathcal{K}]/({v^2_i=u_i v_i=0} ,  1\leqslant i \leqslant m).$$

\noindent Then $R^{*}(\mathcal{K})$ is a bigraded differential algebra with an additive basis
$\{u_J v_I\}$, where $\, I \in \mathcal{K},\,\, J \subseteq [m],\,\, I \cap J = \varnothing;$  $$\mathrm{bideg} \, u_i = (-1, 2), \, \mathrm{bideg} \, v_i = (0, 2),\, d u_i = v_i, \, d v_i = 0.$$

\noindent It is convenient to consider the
$\mathbb{Z} \oplus \mathbb{Z}^m$-grading on $R^*(\mathcal{K})$:  $$\mathrm{mdeg} \, u_i = (-1; 2 e_i), \, \mathrm{mdeg} \, v_i = (0;2 e_i),$$ where $e_i, \, i = 1, \ldots , m,$ are the elements of the standard basis in $\mathbb{Z}^m$.

The multigrading of the algebra of cellular cochains $C^*(\mathcal{Z_K})$ induced from the standard cell decomposition plays an important role in the proof of the existence of a nontrivial Massey product. Consequently, we have the multigrading of the algebra $H^*(\mathcal{Z_K})$.

\medskip

\begin{theorem}[{\cite[Lemma 4.5.3]{1}}]
There is an isomorphism of dg-algebras:
$$R^*(\mathcal{K}) \cong C^*(\mathcal{Z_K}),$$

\noindent where $C^*(\mathcal{Z_K})$ is the algebra of cellular cochains with a natural multiplication inducing the standard product in cohomology.
\end{theorem}

\medskip

Given $J \subseteq [m]$, define the corresponding \textit{full subcomplex} $\mathcal{K}_J$ of a simplicial complex $\mathcal{K}$ as $$\mathcal{K}_J = \{ \, I \in \mathcal{K} \mid I \subseteq J \, \}.$$ For each $\mathcal{K}_J$ we consider the simplicial cochain complex  $(C^*(\mathcal{K}_J), d)$ with the coaugmentation. The group $C^p(\mathcal{K}_J)$ is a free abelian group with a basis of cochains $\chi_L$, where $\chi_L$ is the characteristic function of a simplex $L \in \mathcal{K}_J,\, |L|=p+1$.

\medskip

\begin{theorem}[{\cite[Theorem 3.2.4]{1}}]
There is an isomorphism of cochain complexes
$(C^*(\mathcal{K}_J), d)$
and
$(R^{*-|J|+1, 2J}(\mathcal{K}), d)$:
$$
\xymatrix{
0 \ar[r] & \mathbb{Z} \ar[r]^-{d} \ar[d]_{f_{-1}}^{\cong} &  C^0(\mathcal{K}_J) \ar[r]^-{d} \ar[d]_{f_0}^{\cong} & \cdots \ar[r]^-{d} & C^{p-1}(\mathcal{K}_J) \ar[r]^-{d} \ar[d]_{f_{p-1}}^{\cong} & \cdots \\
0 \ar[r] & R^{-|J|, \, 2J} \ar[r]^-{d} &  R^{-|J|+1, \, 2J} \ar[r]^-{d} & \cdots \ar[r]^-{d} & R^{-|J|+p, \, 2J} \ar[r]^-{d} & \cdots
}
$$
\noindent where $f_p(\chi_L)=\varepsilon(L, J) u_{J \setminus L} v_L$, $\varepsilon(L, J) = \pm 1$  is a certain sign.
\end{theorem}

\medskip

In this way we have an isomorphism of differential graded algebras $$C^*(\mathcal{Z_K}) \cong R^*(\mathcal{K}) \cong \bigoplus \limits_{p \geqslant 0, \, J \subseteq [m]} C^{p-1}(\mathcal{K}_J),  \eqno (2.1)$$
and also
$$H^*(\mathcal{Z_K}) \cong H(R^*(\mathcal{K})) \hm \cong \bigoplus \limits_{p \geqslant 0, \, J \subseteq [m]} \widetilde H^{p-1}(\mathcal{K}_J).$$

The product in the direct sum of simplicial cochains is induced from $R^*(\mathcal{K})$ by the isomorphism (2.7).

\noindent Since the algebras $H^*(\mathcal{Z_K})$ and $H(R^*(\mathcal{K}))$ are multigraded we have:

$$H^{-p, 2J}(\mathcal{Z_K}) \cong H^{-p, 2J}(R^*(\mathcal{K})) \cong \widetilde H^{-p+|J|-1}(\mathcal{K}_J) \subset H^{-p + 2|J|}(\mathcal{Z_K}).$$

\medskip

\begin{theorem}[{\cite[Proposition 3.2.10]{1}}]

The product in $\bigoplus \limits_{p \geqslant 0, \, J \subseteq [m]} C^{p-1}(\mathcal{K}_J)$, induced from $R^*(\mathcal{K})$,
coincides up to a sign with the product defined by the maps

$$
\mu \colon C^{p-1}(\mathcal{K}_I) \times C^{q-1}(\mathcal{K}_J) \to C^{p+q-1}(\mathcal{K}_{I \cup J}),
$$
$$
(\chi_L, \chi_M) \mapsto
\begin{cases}
\chi_{L \sqcup M},&\textit{if  $I \cap J = \varnothing,  L \sqcup M \in \mathcal{K}_{I \sqcup J}$}\\
0,&\textit{otherwise,}\\
\end{cases}, \eqno (2.2)
$$

\noindent where $\chi_L$ is the characteristic function of the simplex $L$.
\end{theorem}

\medskip

A simple $n$-dimensional polytope $P$ is called a \textit{flag polytope} if any set of pairwise intersecting
facets $F_{i_1}, \ldots, F_{i_k}, \, F_{i_s} \cap F_{i_t}  \neq \varnothing$, $s, t = 1, \ldots, k,$ has a nonempty intersection $F_{i_1} \cap \cdots \cap F_{i_k}\neq \varnothing$.

Let $P$ be a simple 3-polytope. Let $F_1, \ldots, F_m$ be its facets. Define a \textit{k-belt} as a cyclic sequence $(F_{i_1}, \ldots, F_{i_k})$ of facets in which only consecutive facets have a nonempty intersection. More precisely: $F_{i_{j_1}} \cap \cdots \cap F_{i_{j_r}} \ne \varnothing$ if and only if $\{\, j_1, \ldots, j_r \, \} \in \{\, \{\, 1, \, 2 \,\}, \, \{\, 2, \, 3 \,\}, \, \ldots, \, \{\, k -1,  \, k \,\}, \, \{\, k, \, 1 \,\} \,\}$. Note that a $k$-belt corresponds to a chordless cycle in the dual complex $\mathcal{K}_P$, which is a triangulation of 2-sphere.

\medskip

\begin{proposition}[{\cite[Proposition 2.3]{8}}]
A simple $3$-polytope $P$ is a flag polytope if and only if  $P \neq \Delta^3$ and $P$ does not contain $3$-belts.
\end{proposition}

\bigskip

\begin{proposition}[{\cite[Proposition 2.5]{8}}]
A simple $3$-polytope P is a flag polytope if and only if any its facet is surrounded
by a $k$-belt, where $k$ is the number of edges of this facet.
\end{proposition}

\bigskip

We say that a polytope $P$ belongs to the \textit{Pogorelov class $\mathcal P$} (or $P$ is a \textit{Pogorelov polytope}) if $P$ is a simple flag 3-dimensional polytope without $4$-belts.
In dimension 3 the class of combinatorial polytopes which admit a right-angled realization in Lobachevsky space $\mathbb{L}^3$ coincides with the Pogorelov class.
Proposition~2.5 implies that $P \in \mathcal P$ if and only if $P \neq \Delta^3$ and $P$ is a simple $3$-polytope without $3$-belts and $4$-belts.

\medskip

\begin{corollary}
A polytope $P$ from the Pogorelov class does not have $3$-gonal and $4$-gonal facets.
\end{corollary}

\medskip

We have the following characteristic property of Pogorelov polytopes.

\medskip

\begin{theorem}[{\cite[Proposition B.2 (b)]{6}}]
A simple $3$-polytope P is a Pogorelov polytope if and only if each pair of its adjacent facets is surrounded by a $k$-belt; if the facets have $k_1$ and $k_2$ edges,
then $k = k_1 + k_2 - 4$.
\end{theorem}

\bigskip
\section{Massey products and Pogorelov polytopes}
\medskip

For moment-angle complex $\mathcal{Z_K}$, a triple Massey product of minimal dimension is given by:

$$H^3(\mathcal{Z_K}) \otimes H^3(\mathcal{Z_K}) \otimes H^3(\mathcal{Z_K}) \hm \to H^8(\mathcal{Z_K}).  \eqno (3.1)$$

These triple Massey products of $3$-dimensional cohomology classes
were completely described by Denham and Suciu in \cite{9}:

\medskip

\begin{theorem}[{\cite[Theorem 6.1.1]{9}}] The following are equivalent:

(1) there exist cohomology classes $\alpha, \beta, \gamma \in H^3(\mathcal{Z_K})$ for which the Massey product $\langle\alpha, \beta, \gamma\rangle$ is defined and nontrivial;

(2) the graph $\mathcal{K}^1$ (the $1$-dimensional skeleton of $\mathcal{K}$) contains an induced subgraph isomorphic to one of the five graphs in Figure 1.
\end{theorem}

\smallskip

\begin{picture}(300,70)
\qbezier(10,11)(35,10)(50,30)
\qbezier(10,11)(15,40)(35,60)
\put(10,11){\circle*{5}}
\put(35, 30){\circle*{5}}
\put(35,45){\circle*{5}}
\put(50,30){\circle*{5}}
\put(50,45){\circle*{5}}
\put(35,60){\circle*{5}}
\put(37,45){\line (1, 0){12}}
\put(37,30){\line (1, 0){12}}
\put(50,32){\line (0, 1){12}}
\put(35,32){\line (0, 1){12}}
\put(35,47){\line (0, 1){12}}
\put(10,10){\line (5, 4){25}}
\put(10,9){\line (2, 3){25}}
\put(50, 45){\line (-1, 1){15}}
\qbezier(71,11)(95,10)(110,30)
\qbezier(71,11)(75,40)(95,60)
\put(71,11){\circle*{5}}
\put(95, 30){\circle*{5}}
\put(95,45){\circle*{5}}
\put(110,30){\circle*{5}}
\put(110,45){\circle*{5}}
\put(95,60){\circle*{5}}
\put(97,45){\line (1, 0){12}}
\put(97,30){\line (1, 0){12}}
\put(110,32){\line (0, 1){12}}
\put(95,32){\line (0, 1){12}}
\put(95,47){\line (0, 1){12}}
\put(70,10){\line (5, 4){25}}
\put(110, 45){\line (-1, 1){15}}
\qbezier(131,11)(155,10)(170,30)
\qbezier(131,11)(135,40)(155,60)
\put(131,11){\circle*{5}}
\put(155, 30){\circle*{5}}
\put(155,45){\circle*{5}}
\put(170,30){\circle*{5}}
\put(170,45){\circle*{5}}
\put(155,60){\circle*{5}}
\put(157,30){\line (1, 0){12}}
\put(170,32){\line (0, 1){12}}
\put(155,32){\line (0, 1){12}}
\put(155,47){\line (0, 1){12}}
\put(130,10){\line (5, 4){25}}
\put(130,9){\line (2, 3){25}}
\put(170, 45){\line (-1, 1){15}}
\qbezier(191,11)(215,10)(230,30)
\qbezier(191,11)(195,40)(215,60)
\put(191,11){\circle*{5}}
\put(215, 30){\circle*{5}}
\put(215,45){\circle*{5}}
\put(230,30){\circle*{5}}
\put(230,45){\circle*{5}}
\put(215,60){\circle*{5}}
\put(217,30){\line (1, 0){12}}
\put(230,32){\line (0, 1){12}}
\put(215,32){\line (0, 1){12}}
\put(215,47){\line (0, 1){12}}
\put(190,10){\line (5, 4){25}}
\put(230, 45){\line (-1, 1){15}}
\qbezier(251,11)(255,40)(275,60)
\put(251,11){\circle*{5}}
\put(275, 30){\circle*{5}}
\put(275,45){\circle*{5}}
\put(290,30){\circle*{5}}
\put(290,45){\circle*{5}}
\put(275,60){\circle*{5}}
\put(277,30){\line (1, 0){12}}
\put(290,32){\line (0, 1){12}}
\put(275,32){\line (0, 1){12}}
\put(275,47){\line (0, 1){12}}
\put(250,10){\line (5, 4){25}}
\put(290, 45){\line (-1, 1){15}}
\end{picture}

\begin{center}
Figure 1.
\end{center}

\medskip
Now we consider the problem of existence of nontrivial Massey products in $H^*(\mathcal{Z_P})$ for Pogorelov polytopes $P$.
As noted in [6, Proposition 4.8], triple Massey products of $3$-dimensional cohomology classes (3.1) are trivial for simple polytopes $P$ without $4$-belts, in particular, for Pogorelov polytopes.
In this paper we prove the following:

\medskip

\begin{theorem}
For any Pogorelov polytope $P$, there exists a nontrivial triple Massey product $\langle \alpha, \beta, \gamma \rangle \subset H^{n+4} (\mathcal{Z_P})$ for some $n \geqslant 5$, where $\alpha \in H^4 (\mathcal{Z_P})$, $\beta \in H^{n-2} (\mathcal{Z_P})$, $\gamma \in H^3 (\mathcal{Z_P})$.
The number $n$ is described in the following lemma.
The indeterminacy of this Massey product is $\alpha \cdot H^n(\mathcal{Z_P})
+ \gamma \cdot H^{n+1}(\mathcal{Z_P})$.
\end{theorem}

\medskip

\medskip

\begin{lemma}
For any Pogorelov polytope $P$, there is a collection of pairwise different facets $\{F_1, \, \ldots, \, F_{l+n-1}\}$ for some natural $n \geqslant 5$ and $l\geqslant 5$ such that the full subcomplex $\mathcal{K}_{\{1, \ldots , l+n-1\}}$ of the complex $\mathcal{K}_P$ has the form shown in Figure 2.
In other words, there exist a triple of facets $F_1, F_2, F_3$ surrounded by a belt.
\end{lemma}

\medskip
\medskip

\hspace{+18ex}
\includegraphics[scale=0.3]{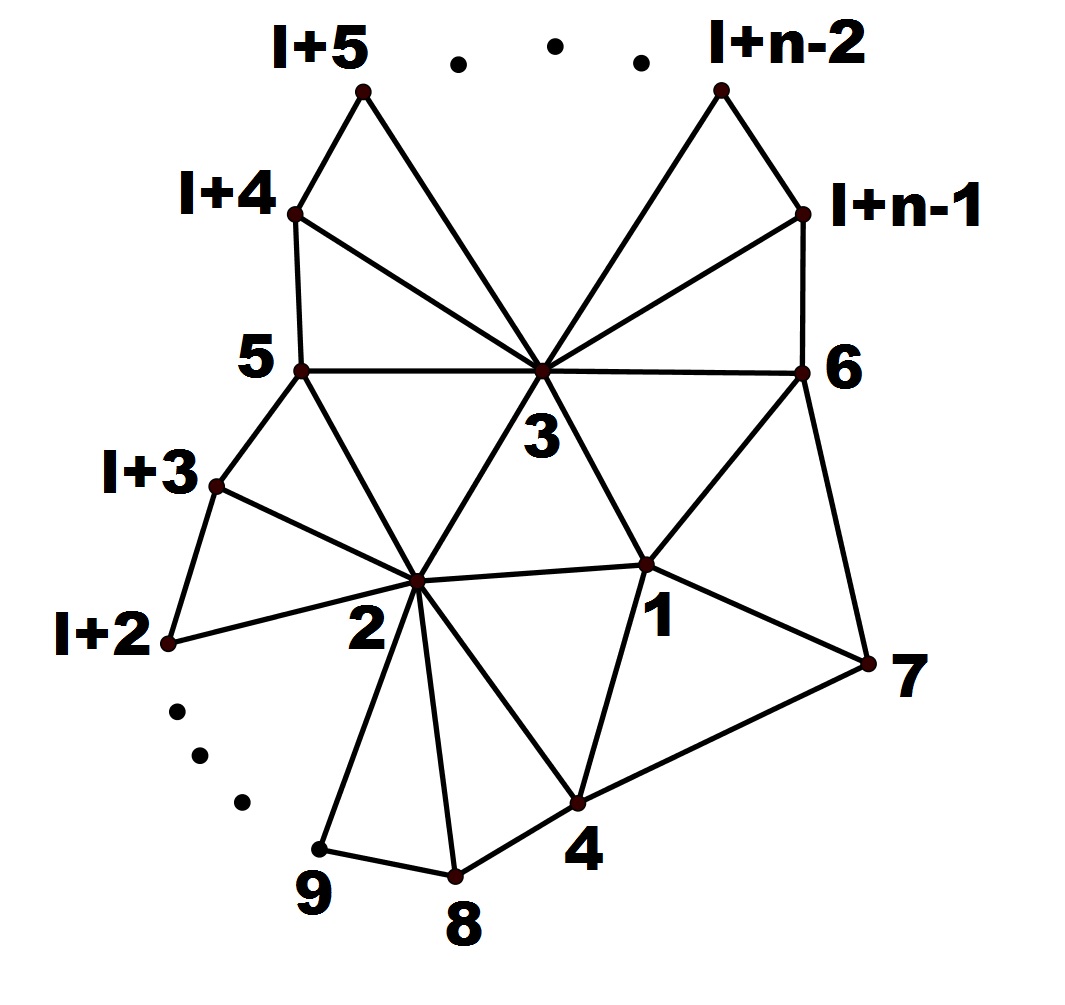}

\begin{center}
Figure 2.
\end{center}

\medskip
\medskip

\begin{proof}
Let $p_k$ be the number of $k$-gonal facets of a simple $3$-polytope $P$. By the Euler theorem, $$3p_3 + 2p_4 + p_5 \hm = 12 + \sum \limits_{k \geqslant 7}(k - 6)p_k.$$ Since $P$ is a Pogorelov polytope, $p_3\hm=0,\, p_4=0$ and hence $p_5 \geqslant 12$. In particular, $p_5 \geqslant 1$.

Choose a $5$-gonal facet $F_1$ of $P$.
Consider an arbitrary vertex $v$ of $F_1$.
Since $P$ is simple, there are exactly 3 facets $F_1, F_2, F_3$ meeting in $v$. Let $F_2$ and $F_3$ be an $l$-gonal facet and an $n$-gonal facet respectively, we denote this by $|F_2|=l,\, |F_3|=n$.
Also, any two facets of $P$ either do not intersect, or intersect at an edge (they are \emph{adjacent facets}).
So $F_1 \cap F_2 = e_{12}$, $F_2 \cap F_3 = e_{23}$, $F_1 \cap F_3 = e_{13}$, $v \in e_{ij}$.
There are two facets intersecting $e_{12} = F_1 \cap F_2 $ at a single vertex. One of them is $F_3$, let $F_4$ be the other one.
The edge $F_2 \cap F_3$ intersects $F_1$ and another facet $F_5$,  and $F_1 \cap F_3$ intersects $F_2$ and $F_6$.
Since $P$ is flag, each facet $F \subset P$ is surrounded by a $k$-belt, where $k=|F|$.
This implies that the facet $F_1$ is surrounded by a $5$-belt and this belt contains $F_2$ and $F_6$.
Therefore, the facets $F_2$ and $F_6$ do not intersect.
Similarly we obtain $F_3 \cap F_4 = \varnothing$, $F_1 \cap F_5 = \varnothing$.
Consider the dual complex $\mathcal{K}_P$, in which the vertex $i$ corresponds to the facet $F_i$.
It follows from the above that the full subcomplex $\mathcal{K}_{\{1, \ldots, 6\}}$ of $\mathcal{K}_P$ has the form shown in Figure 3.

\medskip
\hspace{10ex}
\includegraphics[scale=0.3]{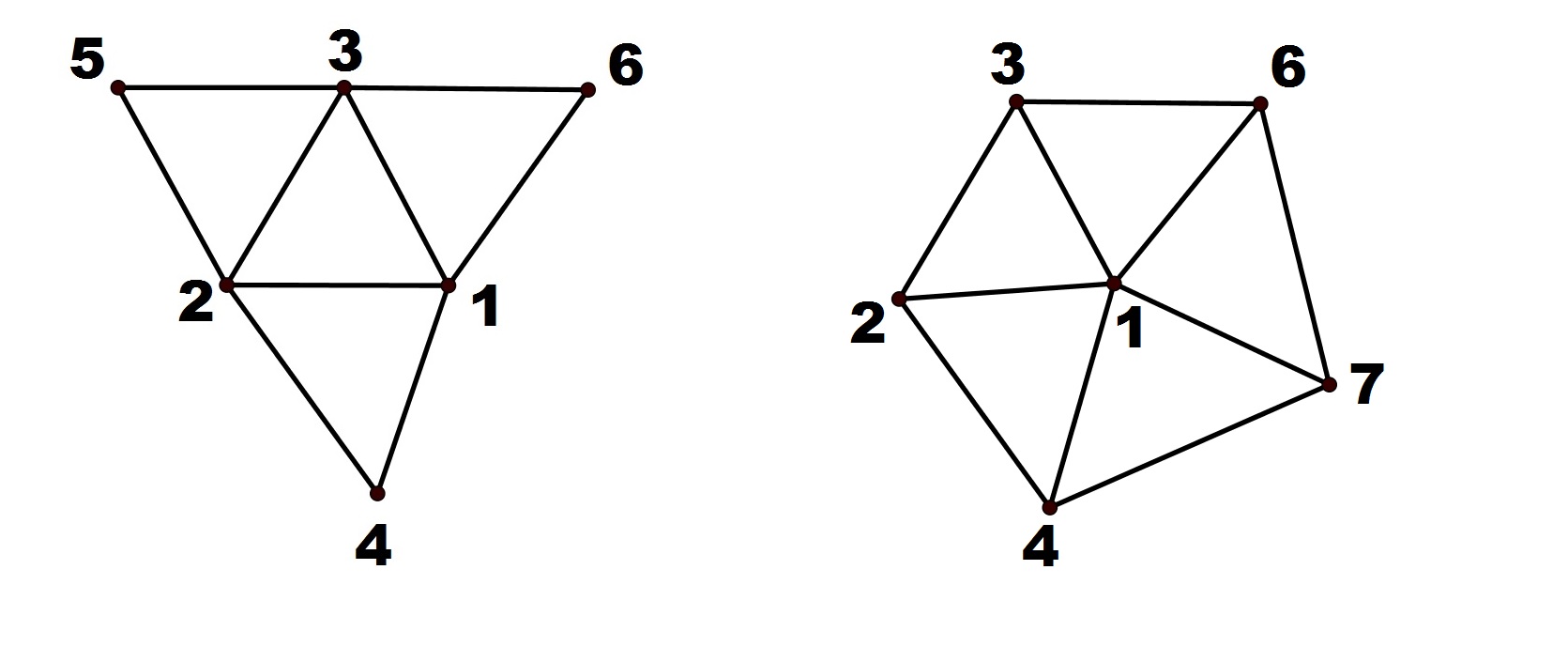}

\begin{center}
\vspace{-2ex}
\hspace{13ex}
Figure 3.
\hspace{25ex}
Figure 4.
\end{center}

\medskip

The pentagonal facet $F_1$ intersects each of the four facets $F_6, F_3, F_2, F_4$ at an edge. The remaining edge of $F_1$ is the intersection of $F_1$ and some other facet $F_7$. The facet $F_1$ is surrounded by a $5$-belt and the simplicial complex $\mathcal{K}_{\{1, 2, 3, 4, 6, 7\}}$ is shown in Figure 4.





Define the set  $\mathcal{G}_2$ consisting of those  facets which intersect $F_2$ and are different from $F_1, \ldots, F_6$ :
$$\mathcal{G}_2=\{F \subset P \mid  F \cap F_2 \neq \varnothing, F\neq F_i, i=1, \ldots, 6 \}.$$
It is easy to see that
$$\mathcal{G}_2=\{F \subset P \mid F \cap F_2 \neq \varnothing, F \cap F_1 = \varnothing,
F \cap F_3 = \varnothing \}.$$
Since $|F_2|=l$ we have $|\mathcal{G}_2|=l-4$. The facets from $\mathcal{G}_2$  are contained in the $l$-belt around $F_2$.
This belt corresponds to a chordless cycle in the dual complex $\mathcal{K}_P$. So we can enumerate the facets from $\mathcal{G}_2$ as follows:
$$\mathcal{G}_2=\{ F_8, \ldots, F_{l+3} \mid F_8 \cap F_4 \neq \varnothing, F_{l+3} \cap F_5 \neq \varnothing \}. $$
The full subcomplex $\mathcal{K}_{\{1, \ldots, 5, 8, \ldots l+3\}}$ is shown in Figure 5.

\medskip

\hspace{7ex}
\includegraphics[scale=0.3]{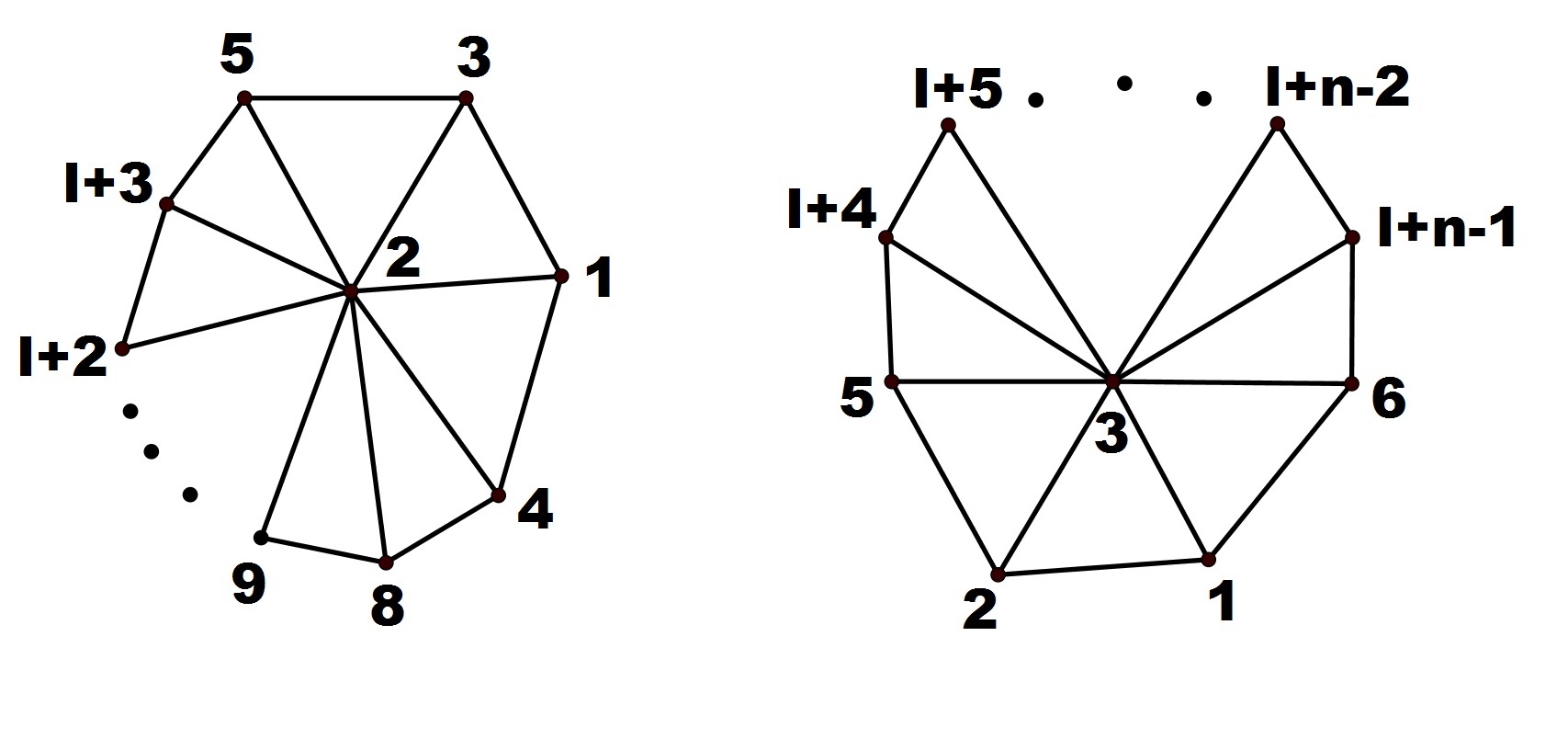}

\vspace{-2ex}
\begin{center}
Figure 5.
\hspace{26ex}
Figure 6.
\end{center}

\medskip
\medskip

 \noindent Similarly, for the facet $F_3$ define the set $\mathcal{G}_3$:
$$\mathcal{G}_3=\{F \subset P \mid  F \cap F_3 \neq \varnothing, F\neq F_i, i=1, \ldots, 6 \}. $$
Since $|\mathcal{G}_3|=n-4$ we also, have that
$$\mathcal{G}_3=\{F \subset P \mid  F \cap F_3 \neq \varnothing, F \cap F_1 = \varnothing,
F \cap F_2 = \varnothing \},$$
$$\mathcal{G}_3=\{ F_{l+4}, \ldots, F_{l+n-1} \mid  F_{l+4} \cap F_5 \neq \varnothing, F_{l+n-1} \cap F_6 \neq \varnothing \}.$$
The full subcomplex $\mathcal{K}_{\{1, 2, 3, 5, 6, l+4, \ldots, l+n-1\}}$ is shown in Figure 6.

Note that the complexes shown in Figures 3, 4, 5 and 6 appear as parts of the complex in Figure 2.
It remains to show that they are patched together in the right way.
That is, we need to show that if $F_i \hm \in \mathcal{G}_2, \, F_j \hm \in  \mathcal{G}_3$ then $F_i \cap F_j = \varnothing, F_7 \cap F_j = \varnothing, F_i \cap F_7 = \varnothing$ ; in particular, $F_i, F_j, F_7$ are different.
Since $P$ is a Pogorelov polytope, Theorem 2.8 implies that the pair of adjacent facets $F_1$ and $F_2$ is surrounded by an $(l+1)$-belt $(F_4, F_8, \ldots, F_{l+3}, F_5, F_3, F_6, F_7)$.
That is, the facets from this sequence are pairwise different and only consecutive facets have a nonempty intersections.
Then since $\mathcal{G}_2= \{ F_8, \ldots, F_{l+3} \}$, we have $F_i \cap F_7 \hm = \varnothing$ if $F_i \in \mathcal{G}_2$.
Similarly considering the pairs of adjacent facets ${F_2, F_3}$ and ${F_3, F_1}$ we
obtain that $\mathcal{G}_2\,\cap\, \mathcal{G}_3 = \varnothing, F_7\, \cap\,F_j = \varnothing$ if $F_j \in \mathcal{G}_3$.
In particular, $(F_4, F_8, \ldots, F_{l+3}, F_5, F_{l+4},$ $ \ldots, F_{l+n-1}, F_6, F_7)$ is a $(l+n-4)$-belt around the triple of facets $\{F_1, F_2, F_3\}$.
\end{proof}

\bigskip

\begin{proof}[Proof of Theorem 3.2]
In the notation of Figure 2,
consider the following three sets of vertices of $\mathcal{K}_P$
(see Figure 7):
$$J_1 = \{5, 6, 7 \}, \, J_2 = \{2, \, l+4, \ldots, \, l+n-1\}, \, J_3 = \{3, 4\}. $$

\medskip

\hspace{-9ex}
\includegraphics[scale=0.3]{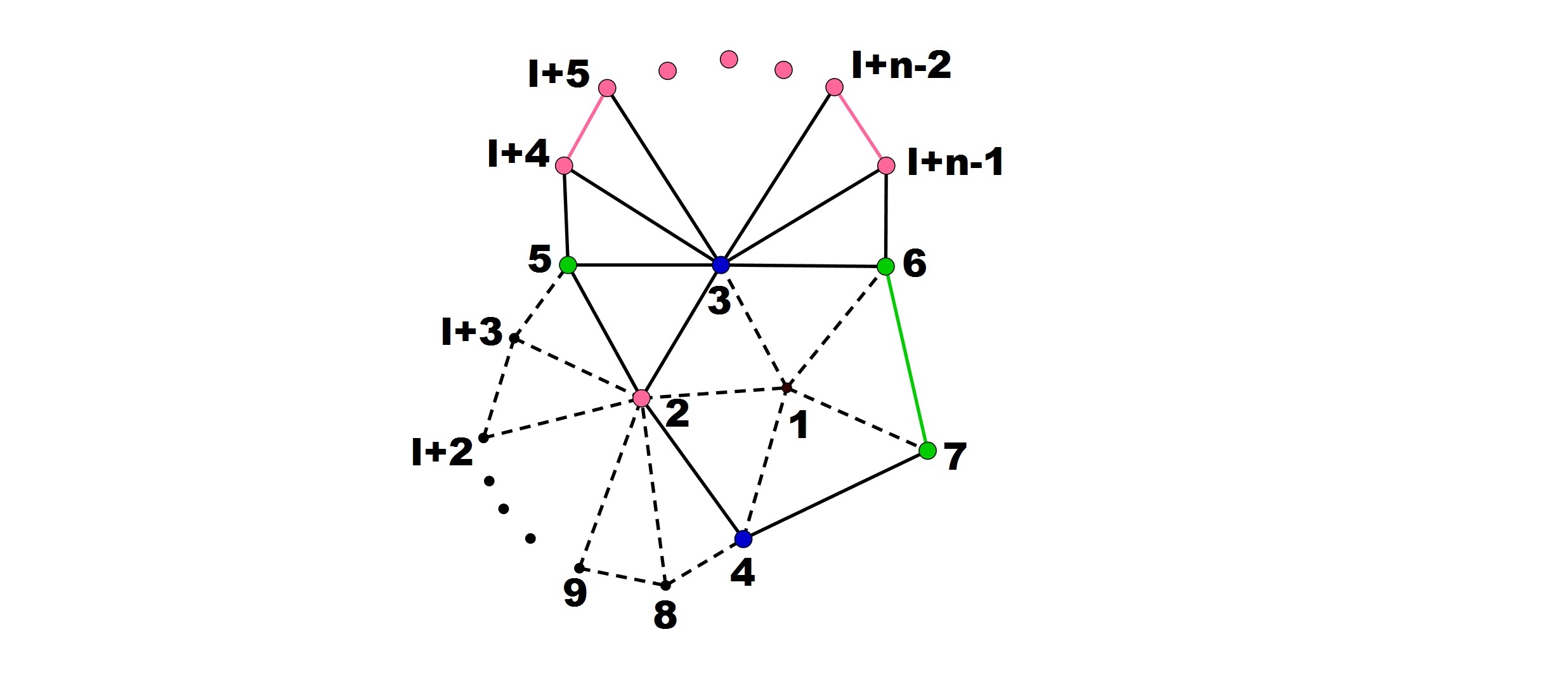}

\vspace{-2ex}
\begin{center}
Figure 7.
\end{center}

\medskip
\medskip

For any $I \in \mathcal{K}_J$, denote by $\chi_{I, J}$ the basis $(|I|-1)$-dimensional simplicial cochain of the complex $\mathcal{K}_J$ taking value $1$ on the simplex $I$.
Now define the following cohomology classes $\alpha, \beta, \gamma$:

$$\alpha = [\chi_{6, J_1} + \chi_{7, J_1}] \in \widetilde H^0(\mathcal{K}_{J_1}) \subset H^4 (\mathcal{Z_P}),$$
$$\beta = [\chi_{2, J_2}] \in \widetilde H^0(\mathcal{K}_{J_2}) \subset H^{n-2} (\mathcal{Z_P}), $$
$$\gamma = [\chi_{4, J_3}] \in \widetilde H^0(\mathcal{K}_{J_3}) \subset H^3 (\mathcal{Z_P}).$$
We consider $\widetilde H^i(\mathcal{K}_{J})$ as subgroups in $H^* (\mathcal{Z_P})$ through the isomorphism (2.7).
Since
$\widetilde H^p(\mathcal{K}_{J})  \cdot \widetilde H^q(\mathcal{K}_{I}) \subset \widetilde H^{p+q+1}(\mathcal{K}_{I \cup J})$,
we obtain that \smallskip
$$\alpha \beta \hm \in \widetilde H^1(\mathcal{K}_{J_1 \cup J_2}), \quad
\beta \gamma \in \widetilde H^1(\mathcal{K}_{J_2 \cup J_3}). $$
We have
$\widetilde H^1(\mathcal{K}_{J_1 \cup J_2}) = \widetilde H^1(\mathcal{K}_{J_2 \cup J_3}) = 0$ because $\mathcal{K}_{J_1 \cup J_2}$ and $\mathcal{K}_{J_2 \cup J_3}$ are contractible.
Thus $\alpha \beta = \beta \gamma = 0$, and the  Massey product $\langle \alpha, \beta, \gamma \rangle$ is defined.
Next,
$$(\chi_{6, J_1} + \chi_{7, J_1}) \cdot \chi_{2, J_2} = 0, \,\,\, \chi_{2, J_2} \cdot \chi_{4, J_3} = \pm \chi_{\{2,4\}, J_2 \cup J_3} = \pm d( \chi_{4, J_2 \cup J_3}), $$
since the product $(2.9)$ in the algebra $\bigoplus \limits_{p \geqslant 0, \, J \subseteq [m]} C^{p-1}(\mathcal{K}_J)$ coincides with the product in $C^*(\mathcal{Z_P})$ up to a sign.
Also, we have
$$(\chi_{6, J_1} + \chi_{7, J_1}) \cdot (\pm \chi_{4, J_2 \cup J_3}) = \pm \chi_{\{4, 7\}, J_1 \cup J_2 \cup J_3}, $$
and
$\pm [\chi_{\{4, 7\}, J_1 \cup J_2 \cup J_3}] \ne 0$ as it is a generator in $H^1(\mathcal{K}_{J_1 \cup J_2 \cup J_3}) =\mathbb{Z}$.
Then
$$\langle \alpha, \beta, \gamma \rangle = \pm [\chi_{\{4, 7\}, J_1 \cup J_2 \cup J_3}] + \alpha \cdot H^n(\mathcal{Z_P})
+ \gamma \cdot H^{n+1}(\mathcal{Z_P})\subset H^{n+4} (\mathcal{Z_P}).$$

We need to prove that the Massey product $\langle \alpha, \beta, \gamma \rangle$ is nontrivial.
Assume the opposite, $0 \in \langle\alpha, \beta, \gamma\rangle$. Then there exist $\nu \in H^n(\mathcal{Z_P})$ and
$\mu \in H^{n+1}(\mathcal{Z_P})$ such that
$$0 \hm = \pm [\chi_{\{4, 7\}, J_1 \cup J_2 \cup J_3}] \, + \, \alpha \cdot \nu \, + \, \gamma \cdot \mu.$$
Since $\,\alpha \in \widetilde H^0(\mathcal{K}_{J_1}),\, \gamma \in \widetilde H^0(\mathcal{K}_{J_3}),$ $ \,
[\chi_{\{4, 7\}, J_1 \cup J_2 \cup J_3}] \in \widetilde H^1(\mathcal{K}_{J_1 \cup J_2 \cup J_3})$,
we can assume that $\nu \in  \widetilde H^0(\mathcal{K}_{J_2 \cup J_3})$,
$\mu \in  \widetilde H^0(\mathcal{K}_{J_1 \cup J_2})$ using the multigrading.
However, $\mathcal{K}_{J_1 \cup J_2}$ and $\mathcal{K}_{J_2 \cup J_3}$ are contractible, hence, $\mu=0, \, \nu=0$.
So $0 = \pm [\chi_{\{4, 7\}, J_1 \cup J_2 \cup J_3}]$, a contradiction.
\end{proof}

\medskip
\medskip

\begin{example}
Let $P$ be a dodecahedron, so $\mathcal{K}_P$ is the boundary of an icosahedron.

\medskip

\medskip
\hspace{20ex}
\includegraphics[scale=0.3]{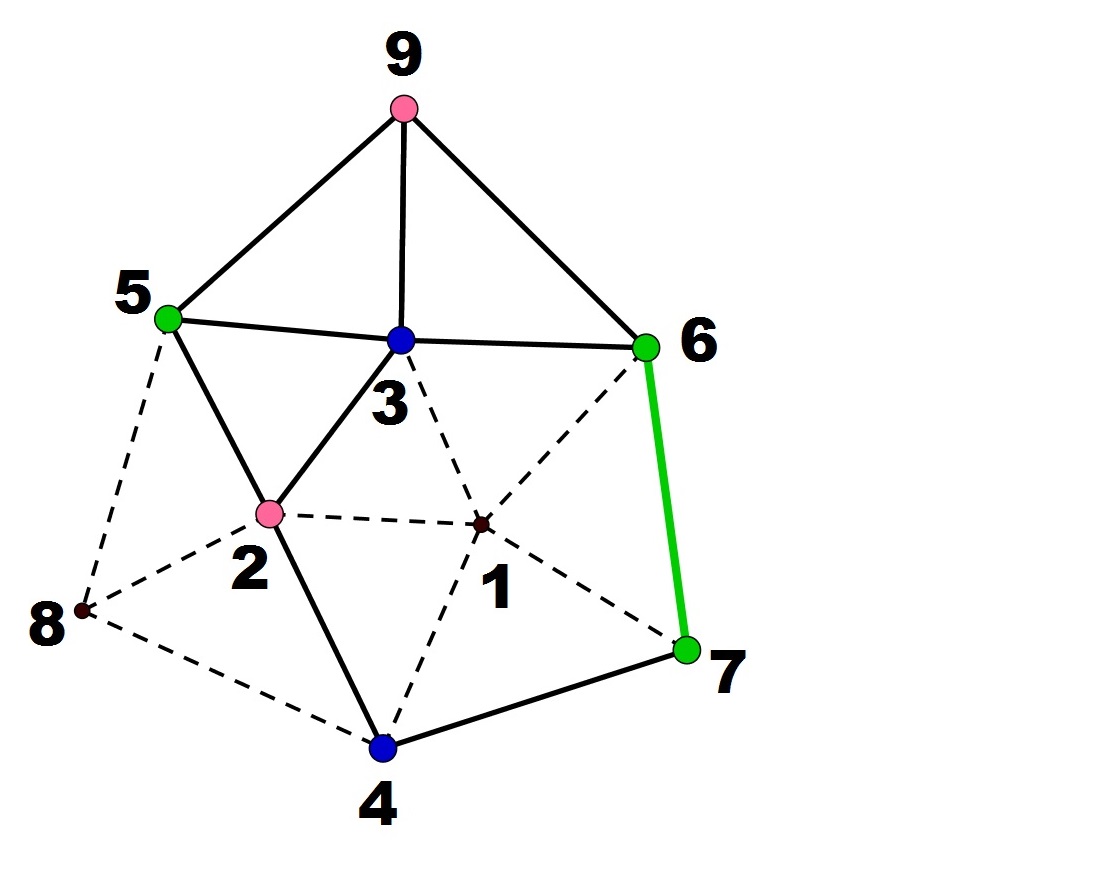}

\vspace{-2ex}
\begin{center}
Figure 8.
\end{center}

\medskip
\medskip

In this case we have
the following sets of vertices of $\mathcal{K}_P$
(see Figure 8):
$$J_1 = \{5, 6, 7 \}, \, J_2 = \{2, 9\}, \, J_3 = \{3, 4\}.$$
The corresponding cohomology classes are
$$\alpha = [\chi_{6, J_1} + \chi_{7, J_1}] \in \widetilde H^0(\mathcal{K}_{J_1}) \subset H^4 (\mathcal{Z_P}),$$
$$\beta = [\chi_{2, J_2}] \in \widetilde H^0(\mathcal{K}_{J_2}) \subset H^{3} (\mathcal{Z_P}),$$
$$\gamma = [\chi_{4, J_3}] \in \widetilde H^0(\mathcal{K}_{J_3}) \subset H^3 (\mathcal{Z_P}).$$
We obtain the following nontrivial Massey product:
$$\langle \alpha, \beta, \gamma \rangle = \pm [\chi_{\{4, 7\}, J_1 \cup J_2 \cup J_3}] \in H^{9} (\mathcal{Z_P}). $$
Note that in the case of dodecahedron we obtain a nontrivial triple Massey product in the lowest possible degree.
\end{example}




\end{document}